\documentclass[11pt,twoside]{article}
\usepackage{amsmath,amssymb,amsthm, txfonts}
 \usepackage{amsfonts}
 \usepackage{amssymb}
 \usepackage{layout}
 \usepackage{verbatim}
 \usepackage{alltt}
 \usepackage{xy}
 \input xy
 \xyoption{all}

\date{}

\newcounter{theorem}[section]

\renewcommand{\thetheorem}{\arabic{section}.\arabic{theorem}}

\newcommand*{\longhookrightarrow}{\ensuremath{\lhook\joinrel\relbar\joinrel\rightarrow}}

\providecommand{\abs}[1]{\lvert#1\rvert}

\newcommand{\C}{\mathbb C}

\newcommand{\Q}{\mathbb Q}
\newcommand{\Z}{\mathbb Z}
\newcommand{\F}{\mathbb F}

\def\p{{\mathfrak p}}
\def\P{{\mathfrak P}}

\def\psl{{\rm PSL}}

\title{On the torsion homology of non-arithmetic hyperbolic tetrahedral groups}  
\author{Mehmet Haluk \c{S}eng\"{u}n}

\begin{document}
\maketitle
\abstract{Numerical data concerning the growth of torsion in the first homology of non-arithmetic hyperbolic tetrahedral groups are collected. 
The data provide support the speculations of Bergeron and Venkatesh on the growth of torsion homology and the regulators for lattices in $\textrm{SL}_2(\C)$.}

\section{Introduction}
In this note, I report on my computations related to the torsion in the first homology of certain non-arithmetic Kleinian groups. 
The motivation for these computations is the recent paper of Bergeron and Venkatesh \cite{bv} on the size of the torsion in the homology of arithmetic uniform lattices in semisimple Lie groups. For $\textrm{SL}_2(\C)$, they have the following result.

\newtheorem*{bv}{Theorem}
\begin{bv} (Bergeron-Venkatesh) Let $\{ \Gamma_n \}$ be a decreasing tower of cocompact arithmetic congruence subgroups of 
$\textrm{SL}_2(\C)$ such that $\bigcap_n \Gamma_n = \{ 1 \}$. Put $X_n= \Gamma_n \backslash \mathbb{H}$ where $\mathbb{H}$ is the hyperbolic 3-space. Then

$$ \lim_{n \rightarrow \infty} \dfrac{\log | H_1(\Gamma_n, E_{k,\ell})_{tor} |}{\textrm{vol}(X_n)} 
                            = \dfrac{1}{6 \pi} \cdot c_{k,\ell}, \ \ \ \ k \not = \ell $$
where $E_{k,\ell}$ is standard module $\textrm{Sym}^k(\Z^2) \otimes \overline{\textrm{Sym}^{\ell}}(\Z^2).$ Here $c_{k,\ell}$ is a positive 
integer depending only on the parameters $k,\ell$, in particular $c_{0,0}=1$.
\end{bv}

The goal of this note is to explore the asymptotic behavior of torsion for non-arithmetic lattices. More precisely, I computationally investigated the limit 
\begin{equation} \lim_{n \rightarrow \infty} \dfrac{\log | H_1(\Gamma_n, \Z)_{tor} |}{\textrm{vol}(X_n)} \end{equation}
for families of groups $\{ \Gamma_n \}$ coming from projective covers of non-arithmetic hyperbolic tetrahedral groups (defined below). 
Note that $\Z \simeq E_{0,0}(\Z)$ in the notation of the theorem. The data I collected show that if one considers only the projective covers  
with vanishing first Betti number (that is, $\textrm{dim}~H_1(\Gamma_n,\Q)= 0 $), then the limit (1) tends to $1/(6 \pi)$. In the case of positive first Betti numbers, the ratios ``log(torsion)/volume" get much smaller than $1/(6 \pi)$. These observations support the general philosophy of Bergeron and Venkatesh as discussed in Sections \ref{sec: bv} and \ref{sec: result} Moreover, the data show that it is extremely rare that residue degree one prime level projective covers of non-arithmetic hyperbolic 3-folds have positive first Betti number. This has been first observed by Calegari and Dunfield in \cite{cd}.

\textbf{Acknowledgements} It is a pleasure to express my gratitude to Grant Lakeland. Upon my request, he kindly computed the matrix realizations which were crucial for my experiments. Most of this work has been done while I was a visitor of the SFB/TR 45 at the Institute for Experimental Mathematics, Essen and I gratefully acknowledge the wonderful hospitality and the state of the art clusters of this institute. I thank Akshay Venkatesh for our encouraging correspondence on the contents of this note. Last but certainly not least, I am grateful to Nicolas Bergeron for hosting me in Paris and explaining to me, among many other things, his joint work with Akshay Venkatesh.

\section{Hyperbolic Tetrahedral Groups} A hyperbolic tetrahedral group is the index two subgroup consisting of orientation-preserving isometries in the discrete group generated by reflections in the faces of a hyperbolic tetrahedron whose dihedral angles are submultiples of $\pi$. Lann\'er \cite{lanner} proved in 1950 that there are thirty two such hyperbolic tetrahedra. 

Consider a tetrahedron in the hyperbolic 3-space $\mathbb{H}$ with vertices $A,B,C,D$ with its dihedral angles submultiples of $\pi$. If the dihedral angles along the edges $AB,AC,BC,DC,DB,DA$ are $\pi / \lambda_1,$ $\pi / \lambda_2,$ $\pi / \lambda_3,$ $\pi / \mu_1,$ $\pi / \mu_2,$ $\pi / \mu_3$ respectively, then we denote the tetrahedron with $T(\lambda_1, \lambda_2, \lambda_3, \mu_1, \mu_2 , \mu_3).$ A presentation for the tetrahedral group $\Gamma$ associated to this tetrahedron is 
$$\Gamma = \langle a,b,c, \mid a^{\mu_1} = b^{\mu_2} = c^{\lambda_3} = (ca)^{\lambda_2}=(cb^{-1})^{\lambda_1}=(ab)^{\mu_3}=1 \rangle.$$

For my computations explicit realizations of these groups in $\psl_2(\C)$ are needed. Note that by Mostow's Rigidity Theorem, 
up to conjugation, such a realization will be unique. A general method to produce a realization is described to me by Grant Lakeland 
in \cite{grant}.

It is well known that there are seven non-arithmetic tetrahedral groups and only of them is cocompact. Moreover, one of the non-cocompact ones is an index two subgroup of another. As nature of the experiment is insensitive to commensurability, I consider only the larger of these two tetrahedral groups. Thus I work with only six groups. Perhaps it should be noted that all the non-cocompact arithmetic tetrahedral groups are commensurable with the Bianchi groups $\textrm{PSL}_2(\Z[\omega])$ with $\omega$ a 4th or 6th root of unity. These two Bianchi groups are covered with the computations in \cite{sen}.

I will now go over each of the six groups, starting with the cocompact one. The volumes of the associated tetrahedra can be found in section 10.4 of the book \cite{egm} or the article \cite{jkrt}. Note that the volume of the 3-fold given by the tetrahedral group is twice the volume of the associated tetrahedron.

\subsection{Descriptions of the Groups}
\label{sec: groups}

\subsubsection{$H(1)$}
Let $H(1)$ be the tetrahedral group attached to the Coxeter symbol
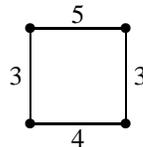
\begin{center}
\setlength{\unitlength}{.5in}
\begin{picture}(2,1.5)
\linethickness{.5 pt}
\put(0,0){\makebox(0,0){$\bullet$}}
\put(1,0){\makebox(0,0){$\bullet$}}
\put(1,1){\makebox(0,0){$\bullet$}}
\put(0,1){\makebox(0,0){$\bullet$}}

\put(0,0){\line(0,1){1}}
\put(0,1){\line(1,0){1}}
\put(1,1){\line(0,-1){1}}
\put(0,0){\line(1,0){1}}

\put(0.50,-0.15){\makebox(0,0)[c]{$4$}}
\put(-0.15,0.50){\makebox(0,0)[c]{$3$}}
\put(0.50,1.15){\makebox(0,0)[c]{$5$}}
\put(1.15,0.50){\makebox(0,0)[c]{$3$}}
\end{picture}   
\end{center}
 The tetrahedron associated to the Coxeter symbol can be described as $T(5,3,2; 4,3,2).$ The volume of the tetrahedron is 
$\simeq 0.3586534401$. A presentation can be given as (see \cite{fgt}) 
$$H(1) = \langle a,b,c \mid a^3=b^2=c^5=(bc^{-1})^3=(ac^{-1})^2=(ab)^4=1 \rangle. $$
where
$$a:=\begin{pmatrix} \frac{2t^3+t^2+t+2}{5} & 1  \\ \frac{-t^3+t^2-2}{5} & \frac{-2t^3-t^2-t+3}{5} \\ \end{pmatrix}, 
\ \ b:=\begin{pmatrix} \frac{-3t^3+t^2-4t+2}{5} & \beta \\ \alpha & \frac{3t^3-t^2+4t-2}{5} \\ \end{pmatrix}, 
\ \  \ \ c:=\begin{pmatrix} t^{-1} & 0  \\ 0 & t \\ \end{pmatrix}$$
where $t$ is a primitive $10$-th root of unity and $\alpha$ is one of the two complex roots of the polynomial
$$x^4+\dfrac{-6t^3+6t^2+8}{5}x^3+\dfrac{-t^3+t^2-3}{5}x^2 + \dfrac{-4t^3+4t^2+2}{25}x + \dfrac{3t^3-3t^2+2}{25}.$$
Moreover, 
$$\beta = (-20t^3+20t^2+35) \alpha^3 + (-50t^3+50t^2+80) \alpha^2 +(9t^3-9t^2-17) \alpha -4t^3+4t^2+6.$$

\subsubsection{$H(2)$}
Let $H(2)$ be the tetrahedral group attached to the Coxeter symbol
\begin{center}
\setlength{\unitlength}{.5in}
\begin{picture}(3,0.5)
\linethickness{.5 pt}
\put(0,0){\makebox(0,0){$\bullet$}}
\put(0.5,0){\makebox(0,0){$\bullet$}}
\put(1,0){\makebox(0,0){$\bullet$}}
\put(1.5,0){\makebox(0,0){$\bullet$}}

\put(0,0){\line(1,0){1.5}}

\put(0.25,0.15){\makebox(0,0)[c]{$6$}}
\put(0.75,0.15){\makebox(0,0)[c]{$3$}}
\put(1.25,0.15){\makebox(0,0)[c]{$5$}}

\end{picture}   
\end{center}
The tetrahedron associated to the Coxeter symbol can be described as $T(5,2,2~;6,2,3).$ It has one ideal vertex. The volume of the tetrahedron is $\simeq 0.1715016613$. A presentation is:
$$H(2) = \langle a,b,c, \mid a^6=b^2=c^2=(ca)^2=(cb^{-1})^5=(ab)^3 =1 \rangle,$$
where
$$       a:=\begin{pmatrix} \zeta & 0  \\ 0 & \zeta^{-1} \\ \end{pmatrix}, 
 \ \     b:=\begin{pmatrix} i & -(\frac{1+\sqrt{5}}{2})i \\ 0 & -i \\ \end{pmatrix}, 
\ \  \ \ c:=\begin{pmatrix} 0 & i  \\ i & 0 \\ \end{pmatrix}.$$

\subsubsection{$H(3)$}
Let $H(3)$ be the tetrahedral group given by the Coxeter symbol 
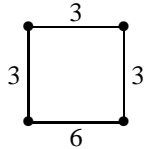
\begin{center}
\setlength{\unitlength}{.5in}
\begin{picture}(2,1.1)
\linethickness{.5 pt}
\put(0,0){\makebox(0,0){$\bullet$}}
\put(1,0){\makebox(0,0){$\bullet$}}
\put(1,1){\makebox(0,0){$\bullet$}}
\put(0,1){\makebox(0,0){$\bullet$}}

\put(0,0){\line(0,1){1}}
\put(0,1){\line(1,0){1}}
\put(1,1){\line(0,-1){1}}
\put(0,0){\line(1,0){1}}

\put(0.50,-0.15){\makebox(0,0)[c]{$6$}}
\put(-0.15,0.50){\makebox(0,0)[c]{$3$}}
\put(0.50,1.15){\makebox(0,0)[c]{$3$}}
\put(1.15,0.50){\makebox(0,0)[c]{$3$}}

\end{picture}   
\end{center}
The tetrahedron associated to the Coxeter symbol can be described as $T(3,3,2~;6,3,2).$ It has two ideal vertices. The volume of the tetrahedron is 
$\simeq 0.3641071004$. A presentation is
$$H(3) = \langle a,b,c, \mid a^6=b^3=c^2=(ca)^3=(cb^{-1})^3=(ab)^2 =1 \rangle,$$
where
$$       a:=\begin{pmatrix} \zeta & i  \\ 0 & \zeta^{-1} \\ \end{pmatrix}, 
 \ \     b:=\begin{pmatrix} \zeta^2 & -i \\ 0 & \zeta^{-2} \\ \end{pmatrix}, 
\ \  \ \ c:=\begin{pmatrix} 0 & i  \\ i & 0 \\ \end{pmatrix}.$$
Here $\zeta = e^{i (\pi/6)}.$

\subsubsection{$H(4)$}
Let $H(4)$ be the tetrahedral group given by the Coxeter symbol 
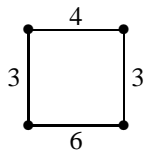
\begin{center}
\setlength{\unitlength}{.5in}
\begin{picture}(2,1.1)
\linethickness{.5 pt}
\put(0,0){\makebox(0,0){$\bullet$}}
\put(1,0){\makebox(0,0){$\bullet$}}
\put(1,1){\makebox(0,0){$\bullet$}}
\put(0,1){\makebox(0,0){$\bullet$}}

\put(0,0){\line(0,1){1}}
\put(0,1){\line(1,0){1}}
\put(1,1){\line(0,-1){1}}
\put(0,0){\line(1,0){1}}

\put(-0.15,0.50){\makebox(0,0)[c]{$3$}}
\put(0.50,1.15){\makebox(0,0)[c]{$4$}}
\put(1.15,0.50){\makebox(0,0)[c]{$3$}}
\put(0.50,-0.15){\makebox(0,0)[c]{$6$}}
\end{picture}   
\end{center}
The tetrahedron associated to the Coxeter symbol can be described as $T(4,3,2~;6,3,2).$ It has two ideal vertices. The volume of the tetrahedron is 
$\simeq 0.5258402692$. A presentation is
$$H(4) = \langle a,b,c, \mid a^6=b^3=c^2=(ca)^3=(cb^{-1})^4=(ab)^2 =1 \rangle,$$
where
$$       a:=\begin{pmatrix} \zeta & i  \\ 0 & \zeta^{-1} \\ \end{pmatrix}, 
 \ \     b:=\begin{pmatrix} \zeta^2 & -\sqrt{-2} \\ 0 & \zeta^{-2} \\ \end{pmatrix}, 
\ \  \ \ c:=\begin{pmatrix} 0 & i  \\ i & 0 \\ \end{pmatrix}.$$
Here $\zeta = e^{i (\pi/6)}.$

\subsubsection{$H(5)$}
Let $H(5)$ be the tetrahedral group given by the Coxeter symbol 
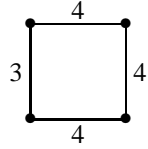
\begin{center}
\setlength{\unitlength}{.5in}
\begin{picture}(2,1.1)
\linethickness{.5 pt}
\put(0,0){\makebox(0,0){$\bullet$}}
\put(1,0){\makebox(0,0){$\bullet$}}
\put(1,1){\makebox(0,0){$\bullet$}}
\put(0,1){\makebox(0,0){$\bullet$}}

\put(0,0){\line(0,1){1}}
\put(0,1){\line(1,0){1}}
\put(1,1){\line(0,-1){1}}
\put(0,0){\line(1,0){1}}

\put(-0.15,0.50){\makebox(0,0)[c]{$3$}}
\put(0.50,1.15){\makebox(0,0)[c]{$4$}}
\put(1.15,0.50){\makebox(0,0)[c]{$4$}}
\put(0.50,-0.15){\makebox(0,0)[c]{$4$}}
\end{picture}   
\end{center}
The tetrahedron associated to the Coxeter symbol can be described as $T(3,4,2~;4,4,2).$ It has two ideal vertices. The volume of the tetrahedron is 
$\simeq 0.5562821156$. A presentation is
$$H(5) = \langle a,b,c, \mid a^4=b^4=c^2=(ca)^4=(cb^{-1})^3=(ab)^2 =1 \rangle,$$
where
$$       a:=\begin{pmatrix} \zeta & \sqrt{-2}  \\ 0 & \zeta^{-1} \\ \end{pmatrix}, 
 \ \     b:=\begin{pmatrix} \zeta & -i \\ 0 & \zeta^{-1} \\ \end{pmatrix}, 
\ \  \ \ c:=\begin{pmatrix} 0 & i  \\ i & 0 \\ \end{pmatrix}.$$
Here $\zeta = e^{i (\pi/8)}.$

\subsubsection{$H(6)$}
Let $H(6)$ be the tetrahedral group given by the Coxeter symbol 
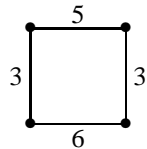
\begin{center}
\setlength{\unitlength}{.5in}
\begin{picture}(2,1.1)
\linethickness{.5 pt}
\put(0,0){\makebox(0,0){$\bullet$}}
\put(1,0){\makebox(0,0){$\bullet$}}
\put(1,1){\makebox(0,0){$\bullet$}}
\put(0,1){\makebox(0,0){$\bullet$}}

\put(0,0){\line(0,1){1}}
\put(0,1){\line(1,0){1}}
\put(1,1){\line(0,-1){1}}
\put(0,0){\line(1,0){1}}

\put(-0.15,0.50){\makebox(0,0)[c]{$3$}}
\put(0.50,1.15){\makebox(0,0)[c]{$5$}}
\put(1.15,0.50){\makebox(0,0)[c]{$3$}}
\put(0.50,-0.15){\makebox(0,0)[c]{$6$}}
\end{picture}   
\end{center}
The tetrahedron associated to the Coxeter symbol can be described as $T(5,3,2~;6,3,2).$ It has two ideal vertices. The volume of the tetrahedron is 
$\simeq 0.6729858045$. A presentation is
$$H(6) = \langle a,b,c, \mid a^6=b^3=c^2=(ca)^3=(c b^{-1})^5=(ab)^2 =1 \rangle,$$
where
$$       a:=\begin{pmatrix} \zeta & i  \\ 0 & \zeta^{-1} \\ \end{pmatrix}, 
 \ \     b:=\begin{pmatrix} \zeta^2 & -(\frac{1+\sqrt{5}}{2})i \\ 0 & \zeta^{-2} \\ \end{pmatrix}, 
\ \  \ \ c:=\begin{pmatrix} 0 & i  \\ i & 0 \\ \end{pmatrix}.$$
Here $\zeta = e^{i (\pi/6)}.$

\section{Projective Covers} \label{sec: proj}
 Let me describe the kind of covers that will be used for the experiment. Let $\Gamma$ be a tetrahedral group as discussed above. As it has finite covolume, there is a number field $K$ with ring of integers $\mathcal{O}$ and a finite set $S$ of finite primes of $\mathcal{O}$ such that there is an embedding 
 $$\Gamma \longhookrightarrow \psl_2(\mathcal{O}_S)$$
 where $\mathcal{O}_S$ is the localization of $\mathcal{O}$ by the primes in $S$. Let $\p$ be a prime ideal of $\mathcal{O}$ not in $S$ and  over the rational prime $p$. Let $\P$ be the prime ideal of $\mathcal{O}_S$ corresponding to $\p$. As the 
ideal is prime, the residue ring $\F$ is a finite field of characteristic $p$. We say that $\p$ is of \textit{residue degree} $a$ if $\F \simeq \F_{p^a}$. Composing the above embedding with the reduction modulo $\P$ map gives us a homomorphism
$$\phi : \ \Gamma \longrightarrow \psl_2(\F).$$
By Strong Approximation Theorem (see Theorem 3.2. of \cite{lackenby}), this map is a surjection for almost all primes $\p$.
Let $B$ denote the Borel subgroup of upper triangular elements in $\psl_2(\F)$. The elements in $\Gamma$ which land in $B$ under $\phi$ 
form a subgroup that we shall denote with $\Gamma_0(\p)$. Note that by rigidity, $\Gamma_0(\p)$ does not depend on $\phi$. 

It is well known that the coset representatives of $\Gamma_0(\p)$ in $\Gamma$ can be identified with the projective line $\mathbb{P}^1(\F)$. 
With this in mind, we will call $\Gamma_0(\p)$ the \textit{projective cover of level $\p$} of $\Gamma$.

\section{More Background for the Experiment}\label{sec: bv}
In this section I will {\it loosely} discuss some aspects of the work of Bergeron and Venkatesh. Any error is due solely to me.

In their recent preprint \cite{bv}, Bergeron and Venkatesh study the growth of torsion in the homology of cocompact arithmetic groups. Their strongest result is in the case where the ambient Lie group is $\textrm{SL}_2(\C)$. In this case one only needs to consider the first homology. Let us specialize our discussion to trivial coefficients (see \cite{bv} p.45). Let $M$ be a hyperbolic 3-manifold of finite volume. The principal object of interest is the product 
$$\abs{H_1(M,\Z)_{tor}} \cdot R(M)$$
where $R(M)$ is the regulator. When $M$ is compact, the regulator $R(M)$ can be given explicitly as 
$\left | \left ( \int_{\gamma_i}\omega_j \right )_{i,j} \right |^{-2}$ where $\omega_1, \hdots, \omega_n$ is an $L^2$-basis of harmonic 1-forms and  $\gamma_1, \hdots , \gamma_n$ is a basis of $H_1(M,\Z)$. In the non-compact case, harmonic 2-forms and $H_2(M,\Z)$ also get involved in $R(M)$.

The general philosophy is that given a tower of finite volume hyperbolic 3-manifolds $M_n$ with ``nice" properties, we should have 
$$\dfrac{\log R(M_n)}{\textit{vol}(M_n)} + \dfrac{\log \abs{H_1(M_n,\Z)_{tor}}}{\textit{vol}(M_n)}  \longrightarrow -\tau^{(2)}(\Z),$$
where $\tau^{(2)}(\Z)$ is a certain quantity which only depends on the ambient Lie group 
and the coefficient module, which is $\Z$ in our limited discussion. This quantity can be explicitly computed and in fact it is equal 
to $-1/(6\pi)$ (see \cite{bv}, p.29). It is a key fact that the product of $\tau^{(2)}$ with $vol(M_n)$ is equal to the $L^2$-analytic torsion of $M_n$.

Hence to study the growth of the torsion, we need to understand that of the regulator. When the first Betti number is zero, the regulator vanishes and according to the general philosophy, the ratio ``log(torsion) /volume" should converge to $1/ (6 \pi)$.  When the first Betti number is positive, the regulator comes into the play and the problem of controlling the conductor 
is a hard one which has connections with the ABC conjecture (see Goldfeld \cite{gold}, p.13). Nevertheless, extensive numerical data 
\cite{sen} show that in the case of arithmetic Kleinian groups, the ratio ``log(torsion)/volume" converges to $1/ (6 \pi)$ regardless the 
first Betti numbers are zero or not. Bergeron and Venkatesh suspect that this has to do with the existence of Hecke operators acting on the homology (see the discussion on p.45 of \cite{bv}). In fact, they conjecture that the ratio ``log(regulator)/volume" goes to zero in the arithmetic case.

\section{The Experiment} 
The method of my computations is essentially the same with the one I employed for the congruence subgroups of the Euclidean Bianchi groups 
in \cite{sen}. Namely, I compute the abelianization of the group instead of its first homology with trivial $\Z$-coefficients. 

Let us fix one of our non-arithmetic tetrahedral groups above, call it $\Gamma$. Let $K$ be the number field generated by the entries of the matrix realizations of the generators of $\Gamma$ given in Section \ref{sec: groups} and let $\mathcal{O}$ be the ring of integers of $K$. 
Let $\p$ be a prime ideal of $\mathcal{O}$ which does not divide the ideal generated by the denominators of the entries of the matrices. 
Then we can reduce the matrices modulo $\p$ as the reductions of the denominators are invertible in the residue field $\F :=\mathcal{O}/ \p$. 
The group $\Gamma_0(\p)$, defined in Section \ref{sec: proj}, has finite index in $\Gamma$ and a set of coset representatives can be identified with the projective line $\mathbb{P}^1(\F)$. The generators of $\Gamma$ act as permutations on the set of coset representatives. The corresponding action on the projective line is the usual action of the matrices, after reducing the entries modulo $\p$. Once the presentation of $\Gamma$ is given to a symbolic algebra program (I used Magma \cite{magma}), next one computes the action of the generator matrices on the projective line. The resulting permutations uniquely determine the subgroup $\Gamma_0(\p)$. Now computing the abelian quotient invariants is a standard functionality.

I have performed two main sets of computations for each of the six tetrahedral groups above. In the first set, I computed the ratios ``log(torsion)/volume" for projective covers with prime level $\p$ over a rational prime $400< p$ with 
norm $N(\p)\leq 50000$. In this way, one captures only prime levels of residue degree one. In a second set, I aimed at prime levels with residue degree at least two. Table \ref{table: ranges} gives the ranges $A,B$ such that $N(\p)\leq A$ and $p \leq B$ for this second set of computations for each group. In addition, I performed an extra set of computations for the groups $H(1)$ and $H(2)$ where 
I considered one prime level over each rational prime $p$ which splits completely in the coefficient field attached to the given matrix realizations of these groups with $p \leq 90000$ and $p \leq 100000$ respectively. 

\begin{table}
\centering
\begin{tabular}{|c|c|c|c|cccc|c|c|c|c|}
\hline
Group & $A$ & $B$ & $\#$ of prime levels    &&&&& Group & $A$ & $B$ & $\#$ of prime levels \\\hline 
$H(1)$& &   &                         &&&&& $H(4)$& 66000 & 257 & 205 \\ \hline
$H(2)$& 150000& 400 &  280            &&&&& $H(5)$& 63000 & 257 & 131 \\ \hline
$H(3)$& 94200 & 307 &  147            &&&&& $H(6)$& 55000 & 233 & 213 \\ \hline
\end{tabular}
\caption{description of the ranges of the second computations} 
\label{table: ranges}
\end{table}

\begin{table}
\centering
\begin{tabular}{|c|c|c|c|c|} \hline
$p$ & $f(\p)$ & $N(\p)$ & rank & log(torsion)/volume  \\ \hline
149& 2& 22201& 35&  0.00223459973429235102460489173120 \\ \hline
157& 2& 24649& 26&  0.00127800797749390070146611206041 \\ \hline
163& 2& 26569& 27&  0.00119388737659470181463861107528 \\ \hline
179& 2& 32041& 42&  0.00284289451885127707642991787490 \\ \hline
179& 2& 32041& 42&  0.00271675896640833678080894982568 \\ \hline
191& 2& 36481& 44&  0.00287523607600174888464895189987 \\ \hline
191& 2& 36481& 44&  0.00298602044559159299824708408103 \\ \hline
193& 2& 37249& 32&  0.00093228187971271537570822454679 \\ \hline
223& 2& 49729& 37&  0.00063399023400905122552042623100 \\ \hline
239& 2& 57121& 56&  0.00204858493468061135230951713718 \\ \hline
239& 2& 57121& 56&  0.00197783049030932370565346415877 \\ \hline
251& 2& 63001& 58&  0.00213916088532177764793377551638 \\ \hline
251& 2& 63001& 58&  0.00207500997945566716013518089664 \\ \hline
269& 2& 72361& 63&  0.00106809775602345757933959528322 \\ \hline
277& 2& 76729& 46&  0.00042737846130476991778286827508 \\ \hline
281& 2& 78961& 65&  0.00115206386564044163781991095396 \\ \hline
283& 2& 80089& 47&  0.00046147256599894360457102615905 \\ \hline
307& 2& 94249& 51&  0.00039717749821458448534705458744 \\ \hline
311& 2& 96721& 72&  0.00161519518746068399522925506475 \\ \hline
311& 2& 96721& 72&  0.00165698128960266506419274897754 \\ \hline
313& 2& 97969& 52&  0.00036262142827441423794453634641 \\ \hline
337& 2& 113569& 56& 0.00031660479790195204470659803953 \\ \hline
359& 2& 128881& 84& 0.00127125170976407173402131001026 \\ \hline
359& 2& 128881& 84& 0.00130261090167121697503469479815 \\ \hline
367& 2& 134689& 61& 0.00025564809254642626286356822328 \\ \hline
373& 2& 139129& 62& 0.00026269400067799622367417655168 \\ \hline
383& 2& 146689& 0&  0.05318174735691597383859191212400  \\ \hline
\end{tabular}
\caption{sample data on growth of torsion for projective covers of H(2) } 
\label{table: homH2}
\end{table}

\section{Results of the Experiment and Observations}\label{sec: result}
The data collected by the experiments are provided on my website \cite{web}. Table \ref{table: homH2} gives sample data for projective covers 
of prime level $\p$ of the tetrahedral group $H(2)$. Here $p$ is the rational prime that $\p$ is over, $f(\p)$ is the residue degree of $\p$,  $N(\p)$ is the norm of $\p$ (which is equal to $p^{f(\p)}$) and rank denotes the first Betti number of the associated 3-fold $M$
(that is, the dimension of $H_1(M,\Q)$).

It is interesting that $\% 99.9$ of the time the homology of a projective cover of residue degree one prime level has first Betti number equal to zero. For each of the six tetrahedral groups, I have computed approximately 5000 projective covers whose levels are prime with residue degree one. For each group, there were no more than 12 such covers with positive first Betti number, all of them with norm less than 200. Note that a similar paucity was observed by Calegari ad Dunfield \cite{cd} among the projective covers of twist knot orbifolds.   

An immediate observation is that if we let $\Gamma$ run only through the projective covers of $G$ with first Betti number equal to {\it zero}, then the data strongly suggest that the limit (1) tends to $1/ (6 \pi) \simeq 0.053051647697298$.
Such a growth was observed in \cite{sen} for the Euclidean Bianchi groups also, {\it but} regardless the first Betti number was zero or not. The convergence here seems to be faster than it was in \cite{sen}. 
 
We see that the behaviour is different in the cases of positive first Betti number. Here the ratios ``log(torsion)/volume"  get much smaller than $1/ (6 \pi)$. This suggests, under the general philosophy of Bergeron and Venkatesh, that the contribution of the regulator comes into play and decreases the contribution of the torsion. Remember that we are in the nonarithmetic case; there is no suitable family of operators acting on the homology. This fact and the numerically observed contribution of the regulator are in favor of the validity of the speculations of Bergeron and Venkatesh on the general behaviour of the regulator.

\end{document}